\documentclass[12pt,a4paper]{article}
\usepackage[cp1251]{inputenc}
\usepackage[T2A]{fontenc}
\usepackage[english]{babel}

\usepackage{amsmath,amsfonts,amssymb,mathrsfs,amsopn,indentfirst,amscd}

\usepackage{graphicx}
\usepackage{amsthm}
\usepackage{hyperref}

\usepackage{amssymb}
\usepackage{amsmath}
\newcommand{\ver}{{\rm ver}}
\newcommand{\vo}{{\rm vol}}

\newcommand{\mes}{{\rm mes_n}}

\newtheorem*{corollary*}{Corollary}

\begin{document}

\title{Geometric Estimates in Interpolation \\ by Linear Functions 
on a Euclidean Ball}
\author{Mikhail Nevskii\footnote{Department of Mathematics,
              P.G.~Demidov Yaroslavl State University, Sovetskaya str., 14, Yaroslavl, 150003, Russia 
              orcid.org/0000-0002-6392-7618 
              mnevsk55@yandex.ru} 
}       
\date{May 9, 2019}
\maketitle

\begin{abstract}
Let $B_n$ be the Euclidean
unit ball in ${\mathbb R}^n$ given by the inequality
$\|x\|\leq 1$, $\|x\|:=\left(\sum\limits_{i=1}^n x_i^2\right)^{\frac{1}{2}}$. 
By  
$C(B_n)$ we mean the space of continuous functions
$f:B_n\to{\mathbb R}$ with the norm
$\|f\|_{C(B_n)}:=\max\limits_{x\in B_n}|f(x)|$. 
The symbol $\Pi_1\left({\mathbb R}^n\right)$ denotes the set of polynomials
in $n$ variables of degree $\leq 1$, i.\,e., the set of linear functions  upon
${\mathbb R}^n$. 
Assume 
$x^{(1)}, \ldots, x^{(n+1)}$ are the vertices
of an 
$n$-dimensional nondegenerate simplex  $S\subset B_n$.
The interpolation projector  
 $P:C(B_n)\to \Pi_1({\mathbb R}^n)$ corresponding to  
$S$ is defined by the equalities
$Pf\left(x^{(j)}\right)=
f\left(x^{(j)}\right).$ Denote by $\|P\|_{B_n}$ the norm of $P$ as an
operator from
$C(B_n)$ onto $C(B_n)$.
We describe the approach 
in which $\|P\|_{B_n}$  can be estimated from below via
the volume of $S$. 

\medskip 

\end{abstract}

 \section{Main Definitions}\label{nev_s1}

We always suppose $n\in{\mathbb N}$. 
An element
$x\in{\mathbb R}^n$ will be written in the form
$x=(x_1,\ldots,x_n).$ 
By definition,  $$\|x\|:=\left(\sum\limits_{i=1}^n x_i^2\right)^{\frac{1}{2}}, \quad 
B_n:=\{x\in{\mathbb R}^n: \|x\|\leq 1 \}, \quad  
Q_n:=[0,1]^n.$$
The notation
$L(n)\asymp M(n)$ means that there exist $c_1,c_2>0$
not depending on $n$ such that
$c_1M(n)\leq L(n)\leq c_2 M(n)$. 
By $\Pi_1\left({\mathbb R}^n\right)$ we mean
the set of polynomials in $n$ variables $\leq 1$, i.\,e.,
the set of linear functions upon ${\mathbb R}^n$.

Let $S$ be a nondegenerate simplex in  ${\mathbb R}^n$ with vertices
 $x^{(j)}=\bigl(x_1^{(j)},\ldots,x_n^{(j)}\bigr)$, $1\leq j\leq n+1$.  
Consider {\it the vertex matrix (or the node matrix)}
$${\bf A} :=
\left( \begin{array}{cccc}
x_1^{(1)}&\ldots&x_n^{(1)}&1\\
x_1^{(2)}&\ldots&x_n^{(2)}&1\\
\vdots&\vdots&\vdots&\vdots\\
x_1^{(n+1)}&\ldots&x_n^{(n+1)}&1\\
\end{array}
\right).$$
We have the equality $\vo(S)=\frac{|\det({\bf A}|}{n!}$.
Put ${\bf A}^{-1}$ 
$=(l_{ij})$. Let us define $\lambda_j$ as polynomials from $\Pi_1\left({\mathbb R}^n\right)$
whose coefficients form the columns of ${\bf A}^{-1}$,  
i.\,e.,~
$\lambda_j(x)$ $:=$
$l_{1j}x_1+\ldots+
l_{nj}x_n+l_{n+1,j}.$
We call $\lambda_j$
{\it the basic Lagrange polynomials} corresponding to this simplex. The numbers $\lambda_j(x)$  
are the barycentric coordinates of a point 
$x\in{\mathbb R}^n$ 
with respect to
$S$. 

For a convex body  $\Omega \subset {\mathbb R}^n$,
denote by
$C(\Omega)$ the space of continuous functions
$f:\Omega\to{\mathbb R}$ with the uniform norm
$$\|f\|_{C(\Omega)}:=\max\limits_{x\in \Omega}|f(x)|.$$
We say that an interpolation projector 
 $P:C(\Omega)\to \Pi_1({\mathbb R}^n)$ corresponds to a simplex
$S\subset \Omega$ if interpolation nodes of $P$ 
coincide with vertices $x^{(j)}$
of this simplex.
Projector $P$ is given by the equalities
$Pf\left(x^{(j)}\right)=
f\left(x^{(j)}\right).$
We have the following analogue of the Lagrange interpolation
formula:
\begin{equation}\label{interp_Lagrange_formula}
Pf(x)=\sum\limits_{j=1}^{n+1}
f\left(x^{(j)}\right)\lambda_j(x). 
\end{equation}
Denote by $\|P\|_\Omega$ the norm of $P$ as an operator from $C(\Omega)$
to $C(\Omega)$. It follows from (\ref{interp_Lagrange_formula}) that
$$\|P\|_\Omega=
\max_{x\in\Omega}\sum_{j=1}^{n+1}
|\lambda_j(x)|.$$
If $\Omega$ is a convex polytope in ${\mathbb R}^n$ (e.\,g., $\Omega=Q_n)$,  
this equality is equivalent to the formula
$$\|P\|_\Omega=
\max_{x\in\ver(\Omega)}\sum_{j=1}^{n+1}
|\lambda_j(x)|,$$
where $\ver(\Omega)$ is the set of vertices of $\Omega$.

Let us define  $\theta_n(\Omega)$ as the minimal value of
$\|P\|_\Omega$
under the conditions $x^{(j)}\in \Omega$. 
In the case $\Omega = Q_n$ 
various relations for the numbers 
$\theta_n(\Omega)$, including the equivalence 
$\theta_n(\Omega)\asymp \sqrt{n}$, were obtained
by the author earlier. This results are systematized in 
[1].
Further on
some estimates were improved (see [2], [3], [4] and references in these papers).

This paper deals with the case $\Omega = B_n$.
We describe the approach when the norm of an interpolation
projector $P: C(B_n) \to \Pi_1({\mathbb R}^n)$ can be estimated from below
through the volume of the corresponding simplex. 
The essential feature of this
approach is the application of the classical Legendre polynomials.
We prove that $ \theta_n (B_n) \asymp \sqrt {n} $. In other words,
the interpolation projector
corresponding to a regular simplex inscribed into the boundary sphere
has the norm equivalent to the minimal possible.

\section{Estimation of $\|P\|_{B_n}$ via the Volume of~$S$}
\label{nev_s2}
{\it The standardized Legendre polynomial of degree 
$n$} is the function
$$\chi_n(t):=\frac{1}{2^nn!}\left[ (t^2-1)^n \right] ^{(n)}$$
(the Rodrigues formula). 
For properties of $\chi_n$, see 
[6],
[7].
The Legendre polynomials are orthogonal on
$[-1,1]$
with respect to the weight $w(t)=1.$ 
It is known that $\chi_n(1)=1$;  if $n\geq 1$, then $\chi_n(t)$ 
strictly increases for~$t\geq 1$.
Denote by
$\chi_n^{-1}$ the function inverse to $\chi_n$ 
on the halfline $[1,+\infty)$.

The appearance of Legendre polynomials in our questions is connected 
with such their property.
For $\gamma\geq 1$, consider the set 
$$E_{n,\gamma}:= \Bigl\{ x\in{\mathbb R}^n :
\sum_{j=1}^n\left |x_j\right| +
\Bigl|1- \sum_{j=1}^n x_j\Bigr|  \leq \gamma \Bigr\}.$$
In 2003 the author established that
\begin{equation}\label{mes_E_eq}
\mes(E_{n,\gamma}) 
= \frac{\chi_n(\gamma)}{n!}. 
\end{equation}
(the proof is also given in [1]).
Utilizing this equality he managed to
obtain the lower estimates for proector's norms related to linear
interpolation on the unit cube $Q_n$. The following relations take place:
\begin{equation}\label{theta_n_nu_h_ineq}
\theta_n(Q_n)
\geq
\chi_n^{-1} \left(\frac{1}{\nu_n}\right)=
\chi_n^{-1} \left(\frac{n!}{h_n}\right).
\end{equation}
Here $\nu_n$ is the maximum volume of a simplex contained in
$Q_n$ and
$h_n$ is the maximum value of $0/1$-determinant of order $n$.
Applying the properties of $\chi_n$ the author got
from (\ref{theta_n_nu_h_ineq})  some more visible inequalities, e.\,g.,
$\theta_n(Q_n)>\frac{1}{e}\sqrt{n-1}.$ 
This estimate occured to be sharp
concerning to
$n$ which led to the equivalence $\theta_n(Q_n)\asymp \sqrt{n}.$

Further on we will extend this approach to linear interpolation of functions
given on the unit ball $B_n$. Denote $\varkappa_n: = \vo(B_n)$. By
$\sigma_n $ we mean the volume of a regular simplex inscribed into $B_n$.

\smallskip
{\bf Theorem 1.}
{\it Assume
$P:C(B_n)\to\Pi_1({\mathbb R}^n)$ is an arbitrary
interpolation projector. Then for
the corresponding simplex 
$S\subset B_n$ and the node matrix ${\bf A}$ we have  
\begin{equation}\label{norm_P_vol_S_ineq}
\|P\|_{B_n}
\geq
\chi_n^{-1} 
\left(\frac{\varkappa_n}{\vo(S)}\right)=
\chi_n^{-1} 
\left(\frac{n!\varkappa_n}{|\det({\bf A})|}\right).
\end{equation}
}

\smallskip
{\it Proof.} 
As known, a regular simplex which is inscribed into a ball
has the maximum possible volume among all simplices being contained in this ball. 
Therefore,
$\vo(S)=\frac{|\det({\bf A})}{n!}\leq \sigma_n$.
For each  $i=1,$ $\ldots,$ $n$, let us subtract from the
$i$th row of ${\bf A}$ its $(n+1)$th row. 
Denote by ${\bf B}$ the submatrix of order $n$ which stands in the first
$n$ rows and columns of the result matrix. Then
$$|\det({\bf B})|=|\det({\bf A})|= n!\vo(S)\leq n!\sigma_n.$$ 
In other words,
\begin{equation}\label{det_B_sigma_n_ineq}
\frac{n!\sigma_n}{|\det({\bf B})|}\geq 1. 
\end{equation}

Let $x^{(j)}$ be the vertices and $\lambda_j$ be the basic Lagrange
polynomials of the simplex $S$.
Since $\lambda_1(x),$ $\ldots$, $\lambda_{n+1}(x)$ are the barycentric
coordinates of a point $x$, we have
$$\|P\|_{B_n}=\max_{x\in B_n} \sum_{j=1}^{n+1}|\lambda_j(x)|
=\max
\Bigl\{ \sum_{j=1}^{n+1}|\beta_j|:
 \, \sum_{j=1}^{n+1}\beta_j=1,
\, 
\sum_{j=1}^{n+1}\beta_jx^{(j)}\in B_n
\Bigr\}.$$
Let us replace $\beta_{n+1}$ with the equal value 
$1-\sum\limits_{j=1}^n
\beta_j.$ The condition 
$\sum\limits_{j=1}^{n+1}\beta_jx^{(j)}$ $\in B_n$ is equivalent to
$\sum\limits_{j=1}^{n}\beta_j(x^{(j)}-x^{(n+1)})\in B^\prime:=B_n-x^{(n+1)}.$ 
This means 
\begin{equation}\label{norm_P_barycentric_from1_ton}
\|P\|_{B_n}=\max\Bigl \{ \sum_{j=1}^{n} \Bigl|\beta_j\Bigr| + 
\Bigl|1-\sum_{j=1}^n\beta_j\Bigr|
\Bigr \}.      
\end{equation}
Maximum in (\ref{norm_P_barycentric_from1_ton})
is taken upon 
$\beta_j$ such that
$\sum\limits_{j=1}^n\beta_j(x^{(j)}-x^{(n+1)})
\in B^\prime.$ 

Consider the nondegenerate linear operator
$F:{\mathbb R}^n \to {\mathbb R}^n$ which maps a point
$\beta=(\beta_1,\ldots,\beta_n)$ into a point
$x=(x_1,\ldots,x_n)$ according to the rule 
$$x=F(\beta)
:=\sum_{j=1}^n \beta_j \left(x^{(j)}-x^{(n+1)}\right).$$ 
We have the matrix equality
$F(\beta)=(\beta_1,\ldots,\beta_n) {\bf B},$
where $\bf B$ is the above introduced 
$(n\times n)$-matrix with the elements $b_{ij}=x_j^{(i)}-x_j^{(n+1)}.$
Put 
$$\gamma^*:=
\chi_n^{-1} 
\left(\frac{n!\varkappa_n}{|\det({\bf B})|}\right).$$
Since $\varkappa_n\geq\sigma_n$,  it follows from 
(\ref{det_B_sigma_n_ineq})   that 
$\gamma^*$ is defined correctly. 
Note the equality 
\begin{equation}\label{chi_n_gamma_star}
\chi_n(\gamma^*)=\frac{n!\varkappa_n}{|\det({\bf B})|}.
\end{equation}

Now suppose $1\leq \gamma < \gamma^*$,
$$E_{n,\gamma}= \Bigl\{ \beta=(\beta_1,\ldots,\beta_n)\in{\mathbb R}^n :
\sum_{j=1}^n\left |\beta_j\right| +
\Bigl|1- \sum_{j=1}^n \beta_j\Bigr|  \leq \gamma \Bigr\}.$$
Let us show that $B^\prime \not \subset F(E_{n,\gamma}).$ 
It is sufficient to get  the inequality $\mes(F(E_{n,\gamma}))
<\mes(B^\prime)=\varkappa_n.$ This is really so:
$$\mes(F(E_{n,\gamma})) <
\mes(F(E_{n,\gamma^*}))=|\det {\bf B}|\cdot \mes(E_{n,\gamma^*})=$$
$$=|\det {\bf B}|\cdot
\frac{\chi_n(\gamma^*)}{n!}=\varkappa_n.$$
We have applied (\ref{mes_E_eq}) 
and (\ref{chi_n_gamma_star}).
Thus, for every $\varepsilon > 0$, there exists a point $x^{(\varepsilon)}$
with the properties
$$x^{(\varepsilon)} 
= \sum_{j=1}^n\beta_j^{(\varepsilon)} \bigl(x^{(j)}-x^{(n+1)}\bigr) \in B^\prime, \quad
\sum\limits_{j=1}^{n}\Bigl|\beta_j^{(\varepsilon)}\Bigr|
+\Bigl|1-\sum\limits_{j=1}^n\beta_j^{(\varepsilon)}\Bigr|
\geq \gamma^*-\varepsilon.$$ 
In view of (\ref{norm_P_barycentric_from1_ton})
this gives $\|P\|_{B_n}\geq \gamma^*-\varepsilon.$
Since $\varepsilon>0$ is an arbitrary, we  obtain
$$\|P\| \geq \gamma^*= 
\chi_n^{-1} \left(\frac{n!\varkappa_n}{|\det({\bf B})|}\right)=
\chi_n^{-1} \left(\frac{n!\varkappa_n}{|\det({\bf A})|}\right)=
\chi_n^{-1} 
\left(\frac{\varkappa_n}{\vo(S)}\right).$$
The theorem is proved.
\hfill$\Box$

\smallskip
{\bf Corollary 1.}
{\it For each $n$,
\begin{equation}\label{theta_n_chi_n_kappa_sigma_ineq}
\theta_n(B_n)
\geq
\chi_n^{-1} 
\left(\frac{\varkappa_n}{\sigma_n}\right).
\end{equation}
}

\smallskip
{\it Proof.} 
Let $P: C(B_n) \to \Pi_1({\mathbb R}^n)$ be an arbitrary interpolation projector. Since the volume
of the corresponding simplex is not greater than $\sigma_n$, inequality  (\ref{norm_P_vol_S_ineq}) 
yields
$$
\|P\|_{B_n}
\geq
\chi_n^{-1} 
\left(\frac{\varkappa_n}{\vo(S)}\right)
\geq
\chi_n^{-1} 
\left(\frac{\varkappa_n}{\sigma_n}\right).
$$
This gives (\ref{theta_n_chi_n_kappa_sigma_ineq}).
\hfill$\Box$

\smallskip
It is known that
\begin{equation}\label{kappa_n_sigma_n_whole}
\varkappa_n=\frac{\pi^
{\frac{n}{2}}}
{\Gamma\left(\frac{n}{2}+1\right)},\qquad
\sigma_n=\frac{1}{n!}\sqrt{n+1}\left(\frac{n+1}{n}\right)^{\frac{n}{2}},
\end{equation}
\begin{equation}\label{kappa_n_even_and_odd}
\varkappa_{2k}=\frac{\pi^{k}}{k!},\qquad
\varkappa_{2k+1}=\frac{2^{k+1}\pi^{k}}{(2k+1)!!}=
\frac{2(k!)(4\pi)^k}{(2k+1)!}.
\end{equation}
Therefore, the estimate (\ref{theta_n_chi_n_kappa_sigma_ineq})
can be made more concrete.

\smallskip
{\bf Corollary 2.}
{\it For every $n$, 
\begin{equation}\label{theta_n_chi_n_through_n_ineq}
\theta_n(B_n)
\geq
\chi_n^{-1} 
\left(\frac{\pi^{\frac{n}{2}}n!}{\Gamma\left(\frac{n}{2}+1\right)\sqrt{n+1}
\left(\frac{n+1}{n}\right)^{\frac{n}{2}}}
\right).
\end{equation}
If $n=2k$, then  (\ref{theta_n_chi_n_through_n_ineq}) is equivalent to the inequality
\begin{equation}\label{theta_n_chi_n_through_n_is_2k_ineq}
\theta_{2k}(B_{2k})
\geq
\chi_{2k}^{-1} 
\left(\frac{\pi^{k}(2k)!}{k!\sqrt{2k+1}
\left(\frac{2k+1}{2k}\right)^k}
\right).
\end{equation}
For $n=2k+1$ we have
\begin{equation}\label{theta_n_chi_n_through_n_is_2k_plus_1_ineq}
\theta_{2k+1}(B_{2k+1})
\geq
\chi_{2k+1}^{-1} 
\left(\frac{2(k!)(4\pi)^{k}}{\sqrt{2k+2}
\left(\frac{2k+2}{2k+1}\right)^{\frac{2k+1}{2}}}
\right).
\end{equation}
}

\smallskip
{\it Proof.} It is sufficient to apply (\ref{theta_n_chi_n_kappa_sigma_ineq}), 
(\ref{kappa_n_sigma_n_whole}), and
(\ref{kappa_n_even_and_odd}).
\hfill$\Box$

\smallskip
{\bf Corollary 3.}
{\it 
Suppose $P$ is a minimal interpolation projector,  i.\,e., 
$\|P\|_{B_n}=\theta_n(B_n)$. Then
\begin{equation}\label{vol_S_for_min_P_ineq}
\vo(S)\geq \frac{\varkappa_n}
{\chi_n\left(\sqrt{n+1}\right)}
=\frac{\pi^
{\frac{n}{2}}}
{\Gamma\left(\frac{n}{2}+1\right)\chi_n\left(\sqrt{n+1}\right)}.
\end{equation}
}

\smallskip
{\it Proof.}
It was proved in [5] that
$\theta_n(B_n)\leq\sqrt{n+1}$. Hence, if $P$ is a minimal projector, then 
$$\sqrt{n+1}\geq \|P\|_{B_n}
\geq
\chi_n^{-1} 
\left(\frac{\varkappa_n}{\vo(S)}\right).$$
We have made use of (\ref{norm_P_vol_S_ineq}). 
It remains to compare the boundary values and take into account
(\ref{kappa_n_sigma_n_whole}).
\hfill$\Box$

\smallskip
Since $|\det({\bf A})|=n!\vo(S)$,  relation (\ref{vol_S_for_min_P_ineq}) 
also  implies an estimate for the determinant of the node matrix corresponding
to a minimal projector. Namely,
$|\det({\bf A})|\geq \varepsilon_n,$ where $\varepsilon_n$ 
is $n!$ times more than the right-hand part of
(\ref{vol_S_for_min_P_ineq}). Without giving the value of $\varepsilon_n$ this 
restriction was used in [5].

\section{The Relation $\theta_n(B_n)\asymp\sqrt{n}$}\label{nev_s3}

The Stirling formula  
$n!=\sqrt{2\pi n}\left(\frac{n}{e}\right)^n
e^{\frac{\zeta_n}{12n}}$,
$0<\zeta_n<1,$
yields 
\begin{equation}\label{n_fact_ineqs}
\sqrt{2\pi n}\left(\frac{n}{e}\right)^n<n!<\sqrt{2\pi n}\left(\frac{n}{e}\right)^n
e^{\frac{1}{12n}}.
\end{equation}
Also we will need the following estimates 
which were proved in 
[1, Section\,3.4.2]: 
\begin{equation}\label{chi_2k_2k_plus_1_ineqs}
\chi_{2k}^{-1}(s)>\left( \frac{(k!)^2s}
{ (2k)!}\right)^{\frac{1}{2k}}, \qquad
\chi_{2k+1}^{-1}(s)>\left( \frac{(k+1)!k!s}
{ (2k+1)!}\right)^{\frac{1}{2k+1}}. 
\end{equation}

\smallskip
{\bf Theorem 2.}
{\it There exist a constant
$c>0$ not depending on $n$ such that
\begin{equation}\label{theta_n_B_n_gt_c_sqrt_n}
\theta_n(B_n)>c\sqrt{n}.
\end{equation}
The inequality 
(\ref{theta_n_B_n_gt_c_sqrt_n}) takes place, e.\,g., with the constant
$$c=
\frac{  \sqrt[3]{\pi} }{\sqrt{12e}\cdot\sqrt[6]{3} } =  0.2135...$$
}

\smallskip
{\it Proof.} First let $n=2k$ be even. 
Making use of
(\ref{theta_n_chi_n_through_n_is_2k_ineq}) and (\ref{chi_2k_2k_plus_1_ineqs}),
we get
$$\theta_{2k}(B_{2k})
\geq
\chi_{2k}^{-1} 
\left(\frac{\pi^{k}(2k)!}{k!\sqrt{2k+1}
\left(\frac{2k+1}{2k}\right)^k}
\right)>
\left(\frac{\pi^{k}(2k)!(k!)^2}{k!\sqrt{2k+1}
\left(\frac{2k+1}{2k}\right)^k(2k)!}
\right)^{\frac{1}{2k}}=$$
$$=
\left(\frac{\pi^{k}k!}{\sqrt{2k+1}
\left(\frac{2k+1}{2k}\right)^k}
\right)^{\frac{1}{2k}}.
$$
Let us estimate $k!$ from below applying  
(\ref{n_fact_ineqs}):
$$\theta_{2k}(B_{2k})>
\left(\frac{\pi^{k}\sqrt{2\pi k}\left(\frac{k}{e}\right)^k}{\sqrt{2k+1}
\left(\frac{2k+1}{2k}\right)^k}
\right)^{\frac{1}{2k}}
=\sqrt{\frac{\pi}{e}} \left(\frac{2\pi k}{2k+1}\right)^{\frac{1}{4k}}
\sqrt{\frac{2k}{2k+1}} \cdot \sqrt{k}.$$
Since $n=2k,$ for each even $n$, 
\begin{equation}\label{complex_frac_even_n}
\theta_{n}(B_{n})>
\sqrt{\frac{\pi}{2e}} \left(\frac{\pi n}{n+1}\right)^{\frac{1}{2n}}
\sqrt{\frac{n}{n+1}} \cdot \sqrt{n}>
\sqrt{\frac{\pi}{3e}} \cdot \sqrt{n}.
\end{equation}
We have utilized the inequalities $\frac{\pi n}{n+1}>1$, $\frac{n}{n+1}>\frac{2}{3}$.
Remark that \linebreak $\sqrt{\frac{\pi}{3e}}=0.6206...$

Now suppose $n$ is odd, i.\,e., $n=2k+1$. In this case
(\ref{theta_n_chi_n_through_n_is_2k_plus_1_ineq})
and (\ref{chi_2k_2k_plus_1_ineqs}) 
give
$$
\theta_{2k+1}(B_{2k+1})
\geq
\chi_{2k+1}^{-1} 
\left(\frac{2(k!)(4\pi)^{k}}{\sqrt{2k+2}
\left(\frac{2k+2}{2k+1}\right)^{\frac{2k+1}{2}}}
\right)>$$
$$>\left(\frac{2(k!)(4\pi)^{k}(k+1)!k!}{\sqrt{2k+2}
\left(\frac{2k+2}{2k+1}\right)^{\frac{2k+1}{2}}(2k+1)!}
\right)^{\frac{1}{2k+1}}=
A^{\frac{1}{2k+1}}\cdot
B^{\frac{1}{2k+1}};
$$
$$A:=\frac{(k!)^2(k+1)!}{(2k+1)!}, \quad
B:=\frac{\sqrt{2}(4\pi)^{k}}{\sqrt{k+1}
\left(\frac{2k+2}{2k+1}\right)^{\frac{2k+1}{2}}}.
$$
From (\ref{n_fact_ineqs}), it follows that 
$$A>\frac{  2\pi k \left(\frac{k}{e}\right)^{2k}\cdot\sqrt{2\pi(k+1)}
\left(\frac{k+1}{e}\right)^{k+1}  }
{\sqrt{2\pi (2k+1)} \left(\frac{2k+1}{e}\right)^{2k+1}e^{\frac{1}{12(2k+1)}}  }=$$
$$=2\pi e^{-k-\frac{1}{12(2k+1)}}
k^{2k+1}(k+1)^{k+1+\frac{1}{2}}(2k+1)^{-\frac{1}{2}-(2k+1)},$$ 
$$ A^{\frac{1}{2k+1}}>
(2\pi)^{\frac{1}{2k+1}} e^{-\frac{k}{2k+1}-\frac{1}{12(2k+1)^2}}
k(k+1)^{\frac{k+1+\frac{1}{2}}{2k+1}}(2k+1)^{-\frac{1}{2(2k+1)}-1}.$$
Also 
$$ (2\pi)^{\frac{1}{2k+1}} >1, \quad \frac{k}{2k+1}+\frac{1}{12(2k+1)^2}<\frac{1}{2},$$ 
$$k(2k+1)^{-\frac{1}{2(2k+1)}-1}=\frac{k}{2k+1} (2k+1)^{-\frac{1}{2(2k+1)}}
>\frac{1}{3}
(2k+1)^{-\frac{1}{2(2k+1)}}>\frac{1}{3\sqrt[6]{3}}. $$
Consequently,  
$$ A^{\frac{1}{2k+1}}>
\frac{1}{3\sqrt[6]{3}\sqrt{e}}(k+1)^{\frac{k+1+\frac{1}{2}}{2k+1}}>
\frac{1}{3\sqrt[6]{3}\sqrt{e}}\sqrt{k+1}.$$
Now let us estimate  
$B^{\frac{1}{2k+1}}$:
$$B^{\frac{1}{2k+1}}= \left( \frac{\sqrt{2}(4\pi)^{k}}
{\sqrt{k+1}\left(\frac{2k+2}{2k+1}\right)^{\frac{2k+1}{2}}}
 \right)^{\frac{1}{2k+1}}=$$
$$=(\sqrt{2})^{\frac{1}{2k+1}}(4\pi)^{\frac{k}{2k+1}}\sqrt{\frac{2k+1}{2k+2}}
(k+1)^{-\frac{1}{2(2k+1)}}>$$
$$>\sqrt[3]{4\pi}\cdot \sqrt{\frac{3}{4}}\cdot 2^{-\frac{1}{6}}=
\sqrt{\frac{3}{2}}\cdot
\sqrt[3]{\pi}.$$
Therefore,
$$
\theta_{2k+1}(B_{2k+1})
>
A^{\frac{1}{2k+1}}\cdot
B^{\frac{1}{2k+1}}>\sqrt{\frac{3}{2}}
\cdot \frac{\sqrt[3]{\pi}}{3\sqrt[6]{3}\sqrt{e}}\cdot\sqrt{k+1}
=\frac{  \sqrt[3]{\pi} }{\sqrt{6e}\sqrt[6]{3} }\cdot\sqrt{k+1}.$$
Since $n=2k+1$, we have $\sqrt{k+1}=
\frac{1}{\sqrt{2}}\sqrt{2k+2}>\frac{1}{\sqrt{2}}\sqrt{n}.$ 
For every odd $n$, 
\begin{equation}\label{complex_frac_odd_n}
\theta_{n}(B_{n})>
\frac{  \sqrt[3]{\pi} }{\sqrt{12e}\cdot\sqrt[6]{3} }\cdot\sqrt{n}.
\end{equation}
The constant $0.2135...$
from the right-hand part of (\ref{complex_frac_odd_n}) 
is
less than the constant from the above inequality
(\ref{complex_frac_even_n}) for even $n$. 
Hence, 
(\ref{complex_frac_odd_n}) 
is true for all $n\in {\mathbb N}$. 
This completes the proof.
\hfill$\Box$



\smallskip
{\bf Corollary 4.}
{\it   $\theta_n(B_n)\asymp \sqrt{n}$.
}

\smallskip
{\it Proof.}
In [5] it was proved that $\theta_n(B_n)\leq \sqrt{n+1}.$
Consequently, the lower estimate $\theta_n(B_n)>c\sqrt{n}$ 
is precise with respect to dimension 
$n$.
\hfill$\Box$

\smallskip
{\bf Corollary 5.}
{\it  Assume $P:C(B_n)\to\Pi_1({\mathbb R}^n)$ is the 
interpolation projector whose nodes coincide with vertices of a regular
simplex being inscribed into the boundary sphere $\|x\|=1$.  
When
$\|P\|_{B_n}\asymp \theta_n(B_n)$.
}

\smallskip
{\it Proof.} As it was shown in [5], 
$\sqrt{n}\leq \|P\|_{B_n}\ \leq \sqrt{n+1}$.   It remains to
utilize the previous corollary.
\hfill$\Box$ 

\smallskip
Our results mean that any interpolation projector $P$
corresponding to an inscrided regular simplex 
has the norm equivalent to the minimal possible. 
The equality $\|P\|_{B_n}=\theta_n(B_n)$ remains proved
only for $1\leq n\leq 4$.
 
\bigskip

The author is grateful to A.\,Yu.~Ukhalov for some useful computer calculations.


\newpage
\centerline{\bf\Large References}

\begin{itemize}

\item[1.]
 Nevskii,~M.\,V.,
{\it Geometricheskie ocenki v polinomialnoi interpolyacii} 
(Geometric Estimates in Polynomial
Interpolation), Yaroslavl': Yarosl. Gos. Univ., 2012 (in~Russian).

\item[2.]
Nevskii,~M.\,V., and Ukhalov, A.\,Yu.,
On numerical charasteristics \linebreak of a simplex and their estimates,
{\it Model. Anal. Inform. Sist.}, 2016, vol.~23, no.~5, pp.~603--619 
(in~Russian).
English transl.: {\it Aut.
Control Comp. Sci.}, 2017, vol.~51, no.~7, pp.~757--769.

\item[3.]
Nevskii,~M.\,V., and Ukhalov, A.\,Yu.,
New estimates of numerical values related to a simplex,
{\it Model. Anal. Inform. Sist.}, 2017, vol.~24, no.~1, pp.~94--110
(in~Russian).
English transl.: {\it Aut.
Control Comp. Sci.}, 2017, vol.~51, no.~7, pp.~770--782.

\item[4.]
Nevskii,~M.\,V., and Ukhalov, A.\,Yu.,
On optimal interpolation by linear functions on an $n$-dimensional cube,
{\it Model. Anal. Inform. Sist.}, 2018, vol.~25, no.~3, pp.~291--311
(in~Russian).
English transl.: {\it Aut.
Control Comp. Sci.}, 2018, vol.~52, no.~7, pp.~828--842.

\item[5.]
Nevskii,~M.\,V., and Ukhalov, A.\,Yu.,
Linear interpolation on a Euclidean ball in ${\mathbb R}^n$
(to appear).

\item[6.]
Szego,~G.,
{\it Orthogonal Polynomials}, New York: American Mathematical Society, 1959.

\item[7.]
Suetin,~P.\,K. {\it Klassicheskie ortogonal'nye mnogochleny},
Moscow: Nauka, 1979 (in~Russian).

\end{itemize}
\end{document}